\documentclass[runningheads]{cl2emult}

\usepackage{graphicx} 

\usepackage{amsmath}   
\usepackage{amsfonts}
\usepackage{mathbbol}
\usepackage{url}
\usepackage{booktabs}
\usepackage[algoruled,vlined]{algorithm2e}
\SetKwInOut{Input}{Input}
\SetKwInOut{Output}{Output}

\newcommand\QQ{\mathbb{Q}}
\newcommand\RR{\mathbb{R}}
\newcommand\Sph{\mathbb{S}}
\newcommand\cF{\mathcal{F}}
\newcommand\cT{\mathcal{T}}
\newcommand\aff{\operatorname{aff}}
\newcommand\bnd{\operatorname{bnd}}
\newcommand\conv{\operatorname{conv}}
\newcommand\size{\operatorname{size}}

\newcommand\SetOf[2]{\left\{\left.#1\,\right|\,#2\right\}}

\newtheorem{thm}{Theorem}
\newtheorem{lem}[thm]{Lemma}
\newtheorem{prp}[thm]{Proposition}
\newtheorem{cor}[thm]{Corollary}

\begin{document}
\title*{Beneath-and-Beyond Revisited}
%
%
%
%
%
\author{Michael Joswig\thanks{Supported by Deutsche Forschungsgemeinschaft, Sonderforschungsbereich~288
    ``Differentialgeometrie und Quantenphysik'' and Forschungszentrum ``Mathematik f\"ur Schl\"usseltechnologien.''}}
%
%
%
\institute{Technische Universit\"at Berlin, Institut f\"ur Mathematik, MA~6--2\\
  Stra\ss{}e des 17.\ Juni 136, 10623 Berlin, Germany\\
  \texttt{joswig@math.tu-berlin.de}, \url{http://www.math.tu-berlin.de/~joswig}}

\maketitle              

\begin{abstract}
  It is shown how the Beneath-and-Beyond algorithm can be used to yield another proof of the equivalence of V- and
  H-representations of convex polytopes.  In this sense this paper serves as the sketch of an introduction to polytope
  theory with a focus on algorithmic aspects.  Moreover, computational results are presented to compare
  Beneath-and-Beyond to other convex hull implementations.
\end{abstract}

\section{Introduction}
\label{sec:Introduction}
One of the most prominent algorithmic problems in computational geometry is the convex hull problem: Given a finite set
of points $S\in\RR^d$, compute the facets (for instance, in terms of their normal vectors) of the convex hull $\conv S$.
To the non-expert it may come as a surprise that fundamental questions concerning the complexity of this problem are
still unsettled.  Numerous methods have been invented through the last three decades.  However, the performance of each
known algorithm heavily depends on specific combinatorial properties of the input polytope~$\conv S$.  Our current
knowledge can largely be summarized by three statements: For many algorithms there is a class of polytopes for which the
given algorithm performs well.  For all known algorithms there is a class of polytopes for which the given algorithm
performs badly.  There are classes of polytopes for which all known algorithms perform badly.

For a comprehensive survey on convex hull algorithms we refer the reader to the paper~\cite{877.68119} by Avis, Bremner,
and Seidel.

The Beneath-and-Beyond algorithm is among the most natural methods for convex hull computation.  It is treated
thoroughly in many text books, see, for instance, Gr\"unbaum~\cite[5.2]{Gruenbaum} and
Edelsbrunner~\cite[Section~8.4]{MR89a:68205}.  The purpose of this paper is twofold: Firstly, we want to sketch how this
particular algorithm can be used to give yet another proof of the ``Main theorem for polytopes'', see
Ziegler~\cite[Theorem~1.1]{GMZ}, which says that the convex hull of finitely many points is exactly the same as a
bounded intersection of finitely many affine halfspaces.  Previously known proofs are based on Fourier-Motzkin
elimination, as in Gr\"otschel~\cite[Kapitel~3]{Groetschel:Skript} and Ziegler~\cite{GMZ}, on the simplex method for
linear optimization, see Schrijver~\cite[Chapter~7]{MR88m:90090}, or non-algorithmic, as in
Gr\"unbaum~\cite[3.1.1]{Gruenbaum} or Matou\v{s}ek~\cite[Theorem~5.2.2]{Matousek:DiscreteGeometry}.  Our strategy for a
proof of the ``Main theorem for polytopes'' via the Beneath-and-Beyond algorithm could be phrased as: Prove everything
directly for simplices and then inductively use triangulations to extend the results to arbitrary polytopes.  Secondly,
we give a brief description of the implementation of the Beneath-and-Beyond algorithm in
\texttt{polymake}~\cite{PolymakeCite,DMV:polymake,SoCG01:polymake}.  The paper closes with a survey of computational
results of Beneath-and-Beyond in comparison with Fourier-Motzkin elimination (or, dually, the double description
method), implemented by Fukuda~\cite{cddlib}, and reverse search, implemented by Avis~\cite{lrs}.  This complements the
computational results in~\cite{877.68119}.

I am indebted to Thilo Schr\"oder who helped to obtain the timing data in Section~\ref{sec:empirical}.  Thanks to
Volker Kaibel, Marc E. Pfetsch, and G\"unter M. Ziegler for their comments and many helpful discussions.  And, finally,
thanks to Ewgenij Gawrilow for co-authoring \texttt{polymake} and his unfailing technical advice.  The polytope images
were produced with \texttt{JavaView}~\cite{javaview} and \texttt{polymake}.

\section{Definitions, an Algorithm, and a Classical Theorem}
\label{sec:Definitions}
A subset $S$ of the Euclidean space~$\RR^d$ is \emph{convex} if for any two points $x,y\in S$ the line segment
$[x,y]=\SetOf{\lambda x+(1-\lambda)y}{0\le\lambda\le1}$ between $x$ and~$y$ is contained in~$S$. We define the
\emph{convex hull} of~$S$, which is denoted by $\conv(S)$, as the smallest convex set containing~$S$.  It is easy to see
that, equivalently,
\begin{equation}
\conv(S)=\SetOf{\sum_{i=1}^n\lambda_ix_i}{x_i\in S,\,\lambda_i\ge0,\,\sum_{j=1}^n\lambda_i=1}.\label{eq:conv}
\end{equation}

For the purpose of this paper the key objects are \emph{(convex) polytopes}, that is, the convex hulls of finite point
sets in~$\RR^d$. 

An affine hyperplane~$H$ defines two (closed) affine halfspaces $H^+$ and~$H^-$ with $H^+\cup H^-=\RR^d$ and $H^+\cap
H^-=H$.  Let $S\subset\RR^d$ be any set.  A hyperplane~$H$ with $H^+\cap S\ne\emptyset$ and $H^-\cap S\ne\emptyset$ is
said to \emph{separate}~$S$.  A hyperplane~$H$ \emph{supports} the set~$S$ if $H$ intersects $S$ non-trivially and if,
moreover, either $S\subset H^+$ or $S\subset H^-$.  We will always assume that $S\subset H^+$, that is, if $H$
supports~$S$ we consider $H$ to be \emph{positively oriented} toward~$S$.

The \emph{dimension} of any set~$S\subset\RR^d$ is defined as the dimension of its affine span~$\aff(S)$; it is denoted
by~$\dim S$.

Throughout the following let $X\subset\RR^d$ be finite.  A \emph{proper face} of the polytope $P=\conv(X)$ is the
intersection of $P$ with a supporting hyperplane which does not contain~$P$.  Note that this last condition is
superfluous in the case where $\dim P=d$, that is, $\aff(P)=\RR^d$.  Occasionally, the empty set and the set~$P$ are
also considered as (non-proper) faces of~$P$.  Some types of faces of a polytope will play particularly important roles
in our investigations: the faces of dimensions $0$, $1$, $d-2$, and~$d-1$, are called \emph{vertices}, \emph{edges},
\emph{ridges}, and \emph{facets}, respectively.

A point $x\in S\subseteq\RR^d$ is a \emph{relatively interior point} of~$S$ if there is a small $\dim S$-dimensional
ball around~$x$ which is contained in~$S$.  The \emph{boundary} of~$S$ is the complement of the relative interior of~$S$
in~$S$; it is denoted by~$\bnd(S)$.

For each facet~$F$ of a polytope~$P$ one can choose a \emph{facet defining} affine halfspace $F^+$ whose boundary
hyperplane supports~$P$ in~$F$.  A \emph{complete} set of facet defining halfspaces contains exactly one such halfspace
for each facet.  If $\dim P=d$ then this choice is unique.

In the following we abbreviate ``$d$-dimensional polytope'' and ``$k$-dimen\-sional face'' by ``$d$-polytope'' and
``$k$-face'', respectively.  The number of $k$-faces of a polytope~$P$ is denoted as $f_k(P)$.  The sequence
$(f_0(P),\dots,f_{d-1}(P))$ is called the $f$-vector of the $d$-polytope~$P$.

Consider $n+1$ affinely independent points $x_1,\dots,x_{n+1}\in\RR^d$.  Their convex hull is an \emph{$n$-simplex}.  A
direct computation yields the following result.

\begin{prp}\label{prp:simplex}
  The $k$-faces of any $n$-simplex $\Delta_n=\conv(x_1,\dots,x_{n+1})$ are precisely the convex hulls of all
  $k+1$~element subsets of $\{x_1,\dots,x_{n+1}\}$.  In particular, the proper $k$-faces are $k$-simplices and
  $f_k(\Delta_n)=\binom{n+1}{k+1}$.  Moreover, each $n$-simplex in~$\RR^n$ is the intersection of its $n+1$ facet
  defining halfspaces, and the boundary of a simplex is the union of its facets.
\end{prp}

Observe that $$\SetOf{\lambda\in\RR^n}{\lambda_i\ge 0,\, \sum\lambda_i=1}$$
is an $(n-1)$-dimensional simplex.  Hence
the Equation~(\ref{eq:conv}) implies that each polytope is the linear projection of a high-dimensional simplex.

A \emph{(geometric) simplicial complex} is a finite collection~$\cT$ of simplices in~$\RR^d$ with the following two
properties:
\begin{enumerate}
\item Each proper face of a simplex in~$\cT$ is also contained in~$\cT$.
\item The intersection of any two simplices $\Delta,\Delta'\in\cT$ is a (possibly empty) face of both, $\Delta$
  and~$\Delta'$.
\end{enumerate}
A simplicial complex~$\cT$ is a \emph{triangulation} of a set $S\subset\RR^d$ if the union of all simplices in~$\cT$
is~$S$.  The $k$-dimensional elements of a simplicial complex~$\cT$ again are called \emph{$k$-faces} of~$\cT$.  A
trivial example: Any simplex (together with its collection of faces) is a triangulation of itself.

Throughout the following we will assume that simplicial complexes (and triangulations) are \emph{pure}, that is, all
(with respect to inclusion) maximal faces have the same dimension.

Usually, it is more convenient to assume that a given polytope~$P=\conv(X)$ affinely spans its ambient
space~$\RR^d$.  This is justified --- also in an algorithmic setting --- due to the following reasoning.  By
performing Gaussian elimination we can determine the dimension $\dim P=\dim X$, and we can even select an affine basis
of the span of~$X$.  Moreover, for instance by omitting redundant coordinates, we can project~$X$ affinely isomorphic to
a linear subspace~$L$ of~$\RR^d$ with $\dim L=\dim X$.  A triangulation of the projection of~$P$ can directly be
lifted back to~$P\subset\RR^d$.  Similarly for facet defining halfspaces.
  
\begin{lem}\label{lem:facet-from-triangulation}
  Let $\cT$ be any triangulation of a $d$-polytope~$P\subset\RR^d$.  Then the following holds.
  \begin{enumerate}
  \item For each face~$F$ of~$P$ the set $\cT(F)=\SetOf{\Delta\in\cT}{\Delta\subseteq F}$ is a triangulation of~$F$.
  \item For each $(d-1)$-face~$\Delta$ of~$\cT$ contained in the boundary of~$P$ there is a unique facet of~$P$ which
    contains~$\Delta$.
  \end{enumerate}
\end{lem}

\begin{proof}
  Let $H$ be a hyperplane which supports~$P$ such that the intersection $H\cap P$ is the face~$F$.  We have to show that
  for each point $x\in F$ there is a face $\Delta\in\cT$ with $x\in\Delta$ and $\Delta\subseteq F$.  For each face
  of~$\cT$ the intersection with~$H$ is again a face of~$\cT$.  But, $\cT$ covers~$P$, that is, some face of~$\cT$
  contains~$x$.  Its intersection with~$H$ is the desired face~$\Delta$.  This proves the first statement.
  
  Now let $\Delta\in\cT$ be a $(d-1)$-face in~$\bnd(P)$.  Since $\cT$ is pure, there exists a $d$-face $\Delta'\in\cT$
  with the property that $\Delta$ is a facet of the $d$-simplex~$\Delta'$.  Therefore, there is a hyperplane~$H$ with
  $H^+\supset\Delta'$ and $H\cap\Delta'=\Delta$.  Choose a point~$y$ in the relative interior of~$\Delta$.  Suppose that
  $H$ separates~$P$, that is, there is a point $x\in H^-\cap P$.  Observe that $\Delta''=\conv(\Delta,x)$ is a
  $d$-simplex, which is contained in~$P$ since $P$ is convex.  Now $\Delta'\cup\Delta''$ is a $d$-dimensional ball which
  contains~$y$ in its interior.  This contradicts $y\in\bnd(P)$, and thus $H$ defines a facet of~$P$.
\end{proof}

A direct consequence of the preceding lemma is the correctness of the Algorithm~\ref{alg:extract} below, which computes
the complete set of facet defining halfspaces of a polytope from a given triangulation.

\begin{algorithm}[H]
  \dontprintsemicolon
  \Input{triangulation~$\cT$ of $P=\conv(X)$}
  \Output{complete set of facet defining halfspaces of~$P$}
  \BlankLine

  $\cF\leftarrow\emptyset$\;
  \ForEach{$(d-1)$-face $\Delta\in\cT$}{
    \If{$\aff(\Delta)$ does not separate~$X$}{
        $\cF\leftarrow\cF\cup\{\aff(\Delta)^+\}$\;
      }
    }
    \Return $\cF$\;
  \caption{\label{alg:extract}Extracting the facets of a polytope from a triangulation.}
\end{algorithm}

Now we are ready to prove the main result of this section: Each polytope is a bounded intersection of finitely many
affine halfspaces.  Observe that the actual statement of the theorem is, in fact, much stronger.  This is necessary to
allow for an easy inductive proof.

\begin{thm}\label{thm:main}
  Let $P=\conv(X)$ be a polytope and $\cF$ a complete set of facet defining affine halfspaces.  Then the following
  holds.
  \begin{enumerate}
  \item There is a triangulation~$\cT$ of~$P$ such that the vertices of~$\cT$ are precisely the points in~$X$.
  \item The polytope~$P$ is the intersection of $\aff(P)$ with the intersection of all halfspaces in~$\cF$.
  \item The boundary of~$P$ is the union of its facets.
  \end{enumerate}
\end{thm}

\begin{proof}
  We give a constructive proof.  More precisely, starting from the finite set~$X$, we construct~$\cT$ and~$\cF$ with the
  desired properties.  As pointed out above we can assume that $P$ is full-dimensional.
  
  For the rest of the proof we fix an arbitrary ordering of the set~$X=\{x_1,\dots,x_n\}$ such that the first~$d+1$
  points $x_1,\dots,x_{d+1}$ are affinely independent.  This ordering gives us a sequence of $d$-polytopes
  $P_k=\conv\{x_1,\dots,x_k\}$ for $k\in\{d+1,\dots,n\}$.  We have $P_{k+1}=\conv(P_k,x_{k+1})$ and we proceed by
  induction on~$k$.
  
  Since $x_1,\dots,x_{d+1}$ are affinely independent, their convex hull~$P_{d+1}$ is a simplex.  The vertices
  $x_1,\dots,x_{d+1}$ of~$P_{d+1}$ clearly are the vertices of the trivial triangulation of~$P_{d+1}$ by itself.  The
  facet defining halfspaces~$\cF_{d+1}$ of $P_{d+1}$ are given by Proposition~\ref{prp:simplex}.  And we have that
  $\bigcap\cF_{d+1}=P_{d+1}$ and $P_{d+1}$ is the union of its facets.
  
  For the inductive step suppose that we have a triangulation $\cT_k$ and the facet defining halfspaces~$\cF_k$
  of~$P_k$.  We construct a triangulation $\cT_{k+1}$ by using Algorithm~\ref{alg:extend} below, where we let $P=P_k$
  and $x=x_{k+1}$.  It is clear that this yields a simplicial complex, but we have to prove that $\bigcup\cT_{k+1}$
  is~$P_{k+1}$.  Without loss of generality $x_{k+1}\not\in P_k$, otherwise $P_{k+1}=P_k$ and the algorithm does not
  change the triangulation.
  
  The induction hypothesis provides us with a facet~$F$ of~$P_k$ which is violated by~$x_{k+1}$.  By
  Lemma~\ref{lem:facet-from-triangulation}a we have a $(d-1)$-face~$\Delta$ of~$\cT_k$ which is contained in~$F$, and
  $\conv(\Delta,x_{k+1})\in\cT_{k+1}$.  Further, consider some point~$x\in P_k$.  By compactness, the line
  segment~$[x,x_k]$ meets the boundary of~$P_k$ in a point.  Due to compactness the line segment~$[x,x_{k+1}]$ meets the
  boundary in a (not necessarily unique) point~$y$.  Again by induction, $\bnd(P_k)$ is the union of its facets, so $y$
  is contained in some facet of~$P_k$ and therefore also in some $(d-1)$-face of~$\cT_k$.  We obtain
  $$P_{k+1}=\bigcup_{x\in P_k}\conv\{x,x_{k+1}\}$$
  and hence $\cT_{k+1}$ is a triangulation of~$P_{k+1}$.
  
  Now the Algorithm~\ref{alg:extract} can be used to determine the complete set~$\cF_{k+1}$ of facet defining halfspaces
  of~$P_{k+1}$.  It is clear that $P_{k+1}\subseteq\bigcap\cF_{k+1}$.  For the reverse inclusion suppose that $z$ is a
  point in $\bigcap\cF_{k+1}\setminus P_{k+1}$.  By repeating the same argument as before we obtain a $(d-1)$-face
  $\Delta\in\cT_{k+1}$ in~$\bnd(P_{k+1})$ such that $\aff(\Delta)$ is a hyperplane which intersects~$P_{k+1}$ in a facet
  and which separates $z$ from~$P_{k+1}$.  This yields the desired contradiction.
  
  Making use of both parts of Lemma~\ref{lem:facet-from-triangulation} allows to conclude that $\bnd(P_{k+1})$ is the
  union of the facets of~$P_{k+1}$.  The theorem is proved.
\end{proof}

\begin{algorithm}[H]
  \dontprintsemicolon
  \Input{triangulation~$\cT$ of $d$-polytope~$P$, facet defining halfspaces~$\cF$ of~$P$, point~$x$}
  \Output{triangulation of~$P'=\conv(P,x)$}
  \BlankLine

  $\cT'\leftarrow\cT$\;
  \ForEach{$F^+\in\cF$}{
    \If{$x\not\in F^+$}{
      \ForEach{\textit{$(d-1)$-face $\Delta\in\cT$ with $\Delta\subset\bnd(F^+)$}}{
        $\cT'\leftarrow\cT'\cup\{\textit{faces of}~\conv(\Delta,x)\}$\;
        }
      }
    }
    \Return $\cT'$\;
  \caption{\label{alg:extend}Extending a triangulation.}
\end{algorithm}

As a direct application we obtain Caratheodory's Theorem:

\begin{cor}\label{cor:Caratheodory} Let $X\subset\RR^d$ be a (full-dimensional) finite set of points and $p\in\conv(X)$.
  Then there is an affinely independent subset $X'=\{x_0,\dots,x_d\}\subseteq X$ with $p\in\conv(X')$.
\end{cor}

\begin{proof}
  By Theorem~\ref{thm:main} the polytope $\conv(X)$ admits a triangulation~$\cT$ with vertex set~$X$.  For $X'$ choose
  the set of vertices of any $d$-simplex in~$\cT$ which contains~$p$.
\end{proof}

This is where we end our algorithmically inspired introduction to polytope theory.  Of course, this is not the end of
the story.  The next step, which then paves the way to the rest of the theory, would be to prove the converse of
Theorem~\ref{thm:main}: Each bounded intersection of finitely many affine halfspaces is a polytope.  One classical way
of proving this is via the Separation Theorem, see Matou\v{s}ek~\cite[Theorem~1.2.4]{Matousek:DiscreteGeometry}, and
duality, see~\cite[5.1]{Matousek:DiscreteGeometry}.  Then Theorem~\ref{thm:main} itself can be applied to prove its
converse.  It is interesting that it is also possible to reverse the order in which one proves these results.


\section{Sizes of Triangulations and Algorithm Complexity}\label{sec:sizes}
One of the major open questions in computational geometry is whether there is a \emph{polynomial total time} convex hull
algorithm, that is, an algorithm whose running time is bounded by a polynomial in the number of vertices \emph{and}
facets.  It is unreasonable to hope for an algorithm whose running time depends polynomially on the input size only,
since there are families of polytopes whose numbers of facets grow exponentially with the number of vertices.  For
instance, take the $d$-dimensional cross polytopes $$\conv\{\pm e_1,\dots,\pm e_d\}$$ with $2d$~vertices and
$2^d$~facets.

In the following we relate the complexity of Beneath-and-Beyond to the size of the triangulations produced.  Then we
review some explicit constructions of polytopes with large triangulations.

\subsection{Complexity Analysis}
We want to examine how the Beneath-and-Beyond method fits into the picture: What is its complexity?  Rather than delving
into technicalities we want to exhibit the geometric core of this question.  The very coarse and schematic description
of the algorithm in steps~\ref{alg:extract} and~\ref{alg:extend} overestimates the costs of the Beneath-and-Beyond
method.  For hints to a more practical approach see the next section.

We begin with Algorithm~\ref{alg:extract}, which extracts the facets of a $d$-polytope~$P$ from a triangulation~$\cT$.
We call the number of $d$-faces of~$\cT$ the \emph{size} of~$\cT$.  If $t=\size(\cT)$ then $\cT$ has at most $(d+1)t$
faces of dimension~$d-1$.  A facet normal vector can be computed from $d$ affinely independent points contained in the
facet by Gaussian elimination, which requires $O(d^3)$~steps.  The overall complexity of Algorithm~\ref{alg:extract} is
then bounded by $O(d^4nt)$, where $n$ is the number of vertices of~$P$.  Now consider Algorithm~\ref{alg:extend}: The
input is a $d$-polytope~$P_k$ with $m_k$~facets and $n_k$~vertices, a triangulation~$\cT_k$ with $t_k$ faces of
dimension~$d$, and one extra point~$x_{k+1}$.  The desired output is a triangulation of~$P_{k+1}=\conv(P_k,x_{k+1})$,
which can be computed by evalutaing $O(dm_kt_k)$ scalar products since each $d$-simplex of~$\cT_k$ contains exactly
$d+1$ simplices of dimension~$d-1$; this gives a total of $O(d^2m_kt_k)$ arithmetic operations.  Summing up we obtain an
upper bound for the complexity of a single Beneath-and-Beyond step.

\begin{lem}
  The vertices, the facets and the triangulation~$\cT_{k+1}$ of $P_{k+1}$ can be computed from the triangulation~$\cT_k$
  of~$P_k$ in $O(d^4\max(m_k,n_k)\,t_k)$ steps.
\end{lem}

We have $t_k\le t_{k+1}$ and $m_k\le(d+1)t_k$.  Setting $m=m_n$ and $t=t_n$ we can sum up for all $k\in\{d+1,\dots,n\}$
to obtain the following result.

\begin{prp}\label{prp:upper-bound}
  The overall complexity of the Beneath-and-Beyond algorithm is bounded by $O(d^5nt^2)$.
\end{prp}

Of course, an implementation of the algorithm as sketched in Algorithms~\ref{alg:extract} and~\ref{alg:extend} is far
from optimal, see Section~\ref{sec:implementation} below for a few more details.  Moreover, the analysis given is very
coarse and could be sharpened easily.  However, for the purpose intended it is good enough: Since the size of the final
triangulation clearly is a lower bound, the Beneath-and-Beyond algorithm runs in polynomial total time if and only if
the size of the triangulation constructed is bounded by a polynomial in $m$ and~$n$.

And this could be the end of the story, because it is known that products of simplices form a family of polytopes where
the size of \emph{any} triangulation is super-polynomial in the number of vertices and facets, see
Haiman~\cite{MR92e:52011} or Avis, Bremner, and Seidel~\cite[Lemma~3]{877.68119}.  We conclude that Beneath-and-Beyond
has a worst case super-exponential total running-time.

It is a consequence of the Upper Bound Theorem, see Ziegler~\cite[Section~8.4]{GMZ}, that the parameter~$t$ is bounded
by $O(n^{\lfloor d/2\rfloor})$.  This bound is actually attained, for instance, by the \emph{cyclic polytopes},
see~\cite[Example~0.6]{GMZ}, which arise as convex hulls of finitely many points on the \emph{moment curve}
$$t\mapsto(t,t^2,t^3,\dots,t^d).$$
If the input points are sorted it is possible to avoid looking at all the facets of~$P_k$ in Algorithm~\ref{alg:extend}.
In this way Algorithm~\ref{alg:extend} can be replaced by a method which takes time which is proportional to the number
of facets of~$P_k$ which are not facets of~$P_{k+1}$.  Amortized analysis then shows that the size~$t$ of the final
triangulation~$\cT_n$ only enters linearly into the cost function.  In particular, for fixed dimension~$d$, one obtains
an $O(n\log n + n^{\lfloor (d+1)/2\rfloor})$ algorithm, see Edelsbrunner~\cite[Section~8.4.5]{MR89a:68205}.

There is a fairly general result due to Bremner~\cite{990.25180} who proves that each incremental convex hull algorithm
has a worst-case super-polynomial total running-time, where \emph{incremental} means that the algorithm has to compute
the convex hulls of all intermediate polytopes~$P_k$.  In particular, this also proves that Beneath-and-Beyond is not a
polynomial total time algorithm.

\subsection{Polytopes with Large Triangulations}
In spite of the fact that there is no polynomial total time convex hull algorithm known, some of the known algorithms
have a polynomially bounded running-time on special classes of polytopes.  In particular, the reverse search method by
Avis and Fukuda~\cite{752.68082} runs in $O(dmn)$ time on simplicial polytopes.  Each simplicial polytope has a small
triangulation: Choose a vertex~$v$ and cone (with apex~$v$) over all the facets not passing through~$v$.  Such a
triangulation is extremely small, it is clearly of size $O(m)$.  Since the size of a triangulation is the decisive
factor for Beneath-and-Beyond's time complexity, this raises the question whether Beneath-and-Beyond could at least have
a polynomial total running-time on simplicial polytopes.  We note that the polytopes in Bremner's
construction~\cite{990.25180} are products of simplicial polytopes, that is, they are neither simplicial nor simple.

A \emph{placing triangulation} of a polytope~$P$ is a triangulation produced by the Beneath-and-Beyond algorithm for
some ordering of the vertices.  For a given ordering, the size of the corresponding placing triangulation is related to
the question of how many facets of the intermediate polytope~$P_k$ are violated by the next vertex~$v_{k+1}$.  A placing
triangulation with respect to a given vertex order is the same as a pushing triangulation with respect to the reverse
order, see Lee~\cite{HDCG:Lee}.  In particular, placing triangulations are lexicographic and thus regular.

In the dual setting, that is with the roles of vertices and facets interchanged, Avis, Bremner, and
Seidel~\cite{877.68119} introduced the concept of \emph{dwarfing}: Informally speaking, a separating hyperplane~$H$ is
called \emph{dwarfing} for a polytope~$P$ if all facets of~$P$ are also facets of~$P\cap H^+$, and if very many vertices
of~$P$ are not vertices of~$P\cap H^+$.  Equivalently, the star of the vertex~$v_H$ dual to~$H$ in any placing
triangulation of the polar polytope~$(P\cap H^+)^*$, where $v_H$ comes last, is large.

\begin{figure}
  \begin{center}
    \includegraphics[width=10cm]{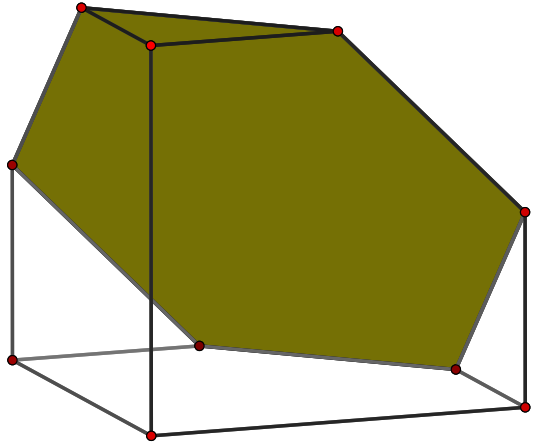}
    \caption{Dwarfed $3$-cube~$c_3$ with dwarfing facet marked.}\label{fig:dwarfed-3cube}
  \end{center}
\end{figure}

\begin{figure}
  \begin{center}
    \includegraphics[width=10cm]{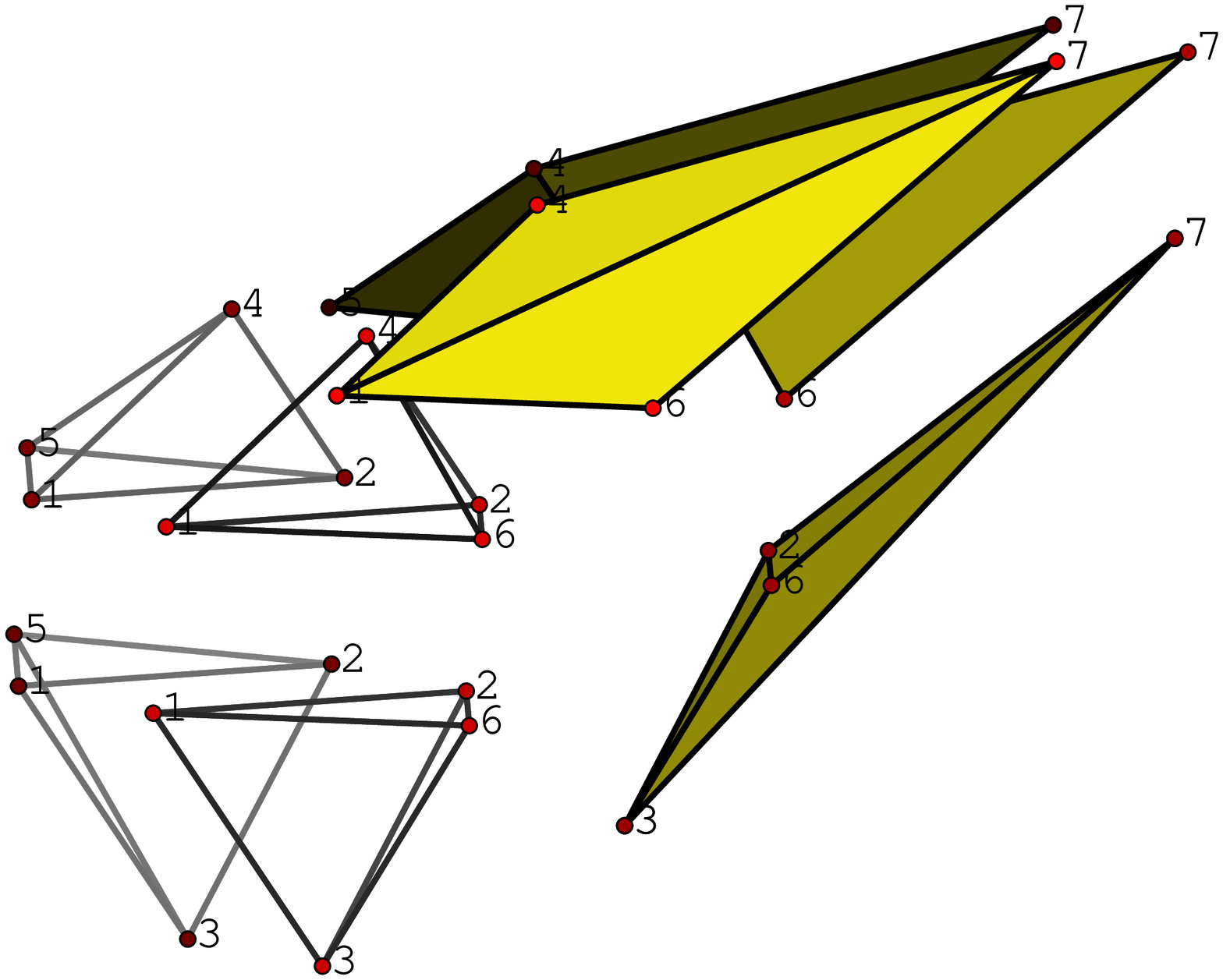}

    \caption{Explosion of the placing triangulation of the polar polytope~$c_3^*$ induced by the ordering indicated as
      vertex labels.  Only the four $3$-simplices containing the final vertex (numbered~$7$), which corresponds to the
      dwarfing facet of~$c_3$, are displayed in solid.}
    \label{fig:polar-dwarfed-3cube}
  \end{center}
\end{figure}

The regular $d$-dimensional cube $C_d=[0,1]^d$ has $2^d$~vertices and $2d$~facets.  Consider the affine halfspace
$H_d^+=\SetOf{x\in\RR^d}{\sum x_i\le3/2}$.  The boundary hyperplane $H_d=\partial H_d^+$ separates the cube~$C_d$.  Now,
as in~\cite[Theorem~4]{877.68119}, call the simple $d$-polytope
$$c_d=C_d\cap H_d^+$$
the \emph{dwarfed $d$-cube}. It has $2d+1$~facets and $d^2+1$~vertices: The origin~$0$ and the
$d$~unit vectors $e_1,\dots,e_d$ are the only vertices of~$C_d$ which are contained in~$H_d^+$; there is one new vertex
$e_i+\frac12e_j$ for each $i,j\in\{1,\dots,d\}$ with $i\ne j$.  For illustrations see Figures~\ref{fig:dwarfed-3cube}
and~\ref{fig:polar-dwarfed-3cube}.  By Barnette's Lower Bound Theorem, see Br{\o}nsted~\cite[\S19]{MR84d:52009}, $d^2+1$
is the minimal number of vertices for a simple $d$-polytope with $2d+1$~facets.  It follows that the dwarfed cubes can
be obtained from a simplex by repeated truncation of faces.

Since the cube and the dwarfed cube have $d+1$ vertices in common, there are $2^d-d-1$ vertices of~$C_d$ which are not
vertices of~$c_d$.  We conclude the following result.

\begin{prp}
  Let $v_1,\dots,v_{2d+1}$ be an ordering of the vertices of the polar $c_d^*$ of the dwarfed $d$-cube~$c_d$ such that
  the last vertex $v'=v_{2d+1}$ corresponds to the dwarfing facet~$H_d$ of~$c_d$.  Then the number of $d$-simplices
  which contain~$v'$ in the induced placing triangulation of~$c_d^*$ equals $2^d-d-1$.
\end{prp}

The dwarfed cubes form a family of polytopes that are ``bad'' as input for the Beneath-and-Beyond algorithm.  Since
there is only one dwarfed cube per dimension, this does not tell anything about the situation in fixed dimension.
However, the same idea can be applied to another (bi-parametric) family of polytopes, whose number of vertices and
facets is unbounded even if the dimension is fixed.  We sketch the construction from~\cite{877.68119}, and we omit the
proofs.

For $d=2\delta\ge 4$ and $s\ge3$ let $G_d(s)$ be the $d$-dimensional polytope defined by the following list of $\delta
s$ linear inequalities, all of which define facets:
\begin{eqnarray}
y_k &\ge& 0\label{eq:product:1}\\
sx_k-y_k &\ge& 0\label{eq:product:2}\\
(2i+1)x_k+y_k &\le& (2i+1)(s+i)-i^2+s^2\label{eq:product:3}\\
(2s-3)x_k+y_k &\le& 2s(2s-3),\label{eq:product:4}
\end{eqnarray}
where $i\in\{0,\dots,s-4\}$, $k\in\{1,\dots,\delta\}$, and a vector in~$\RR^d$ is written as
$(x_1,y_1,x_2,y_2,\dots,x_\delta,y_\delta)$.  The polytope~$G_d(s)$ is the product of $\delta$ copies of the $s$-gon
obtained by fixing~$k$ in the list of inequalities~\ref{eq:product:1} to~\ref{eq:product:4}.  This product of $s$-gons
has $s^\delta$ vertices.  There is a dwarfing halfspace $H_{d,s}^+=\SetOf{x\in\RR^d}{\sum x_i\le 2s-1}$, and the
\emph{dwarfed product of $s$-gons} $$g_d(s)=G_d(s)\cap H_{d,s}^+$$
is simple.  It has $\delta s+1$ facets but only
$(d-1)(\delta s+1-d)+2$ vertices, which, again, meets Barnette's lower bound; see~\cite[Theorem~6]{877.68119}.

Since $G_d(s)$ and $g_d(s)$ can have at most $(d-1)(\delta s+1-d)+1$ vertices in common, this yields a similar result as
for the dwarfed cubes.

\begin{prp}\label{prp:triangulation-dwarfed-polygons}
  Let $v_1,\dots,v_{(d-1)(\delta s+1-d)+2}$ be an ordering of the vertices of the polar $g_d(s)^*$ of the dwarfed
  product of polygons $g_d(s)$ such that the last vertex $v'=v_{(d-1)(\delta s+1-d)+2}$ corresponds to the dwarfing
  facet~$H_{d,s}$ of~$g_d(s)$.  Then the number of $d$-simplices which contain~$v'$ in the induced placing triangulation
  of~$g_d(s)^*$ is at least $s^\delta-(d-1)(\delta s+1-d)-1$.
\end{prp}

This already shows that the worst case running time of Beneath-and-Beyond is not polynomially bounded in combined the
size of the input and the output total time algorithm.  Additionally, one can show that $G_d(s)$ and $g_d(s)$ share
exactly $\delta(s-2)+1$ vertices, so the precise number of $d$-simplices containing~$v'$ is $s^\delta-\delta
s+2\delta-1$.

The closer analysis in~\cite[Theorem~12]{877.68119} reveals that essentially from the preceding
Proposition~\ref{prp:triangulation-dwarfed-polygons} it follows that even a typical placing triangulation (that is, with
respect to a random ordering) of~$g_d(s)^*$ grows super-polynomially with the number of vertices and facets.

\begin{thm}
  The polar dwarfed products of polygons $g_d(s)^*$ are simplicial $d$-polytopes with $(d-1)(ds/2+1-d)+2$~vertices and
  $ds/2+1$~facets, but with an expected size of a placing triangulation of order~$\Omega(s^{d/2}/(d+1))$.
\end{thm}

It seems to be open whether there is a class of simplicial polytopes such that \emph{each} placing triangulation is
large compared to their number of (vertices and) facets.

\section{On the Implementation}\label{sec:implementation}
In the preceding section we related the performance of the Beneath-and-Beyond algorithm to the size of certain
triangulations.  While this captures the main ideas of the method, an implementation which is feasible for practical
problems is slightly more involved.  We sketch the well-known key points below; for a more thorough discussion see
Edelsbrunner~\cite[8.4.5]{MR89a:68205}.

Rather than extracting the facets of the intermediate polytope~$P_k$ from the triangulation~$\cT_k$ from scratch (as
Algorithm~\ref{alg:extract} suggests) it is more natural to store the facets and their neighborhood structure in terms
of the dual graph of~$P_k$.  Adding the next vertex $v_{k+1}$ then requires: (i) to find one violated facet, (ii) to
perform a breadth-first-search in the dual graph of~$P_k$ to find all the other violated facets, (iii) to extend the
triangulation $\cT_k$ of~$P_k$ to the triangulation $\cT_{k+1}$ of~$P_{k+1}$ by coning over the induced triangulation of
the violated facets, (iv) to determine the new facets, that is, the facets of~$P_{k+1}$ which are not facets of~$P_k$
by examining~$\cT_{k+1}$.

Edelsbrunner~\cite[8.4.5]{MR89a:68205} advocates to sort the input points lexicographically in order to eliminate the
time for finding the first violated facet, which then allows to find all the violated facets in linear time (in fixed
dimension).  Our implementation does not rely on ordered input, but it is programmed in a way such that Edelsbrunner's
analysis applies, if the input happens to be ordered.  This more flexible algorithm has the advantage that it can
possibly benefit from smaller placing triangulations.  The best strategy for most cases seems to be to permute the input
randomly rather than to sort it.

We omit a discussion of the data structures required in the implementation.  Instead we want to to spend a few words on
the arithmetic to be used.  In a principal way, convex hull computations makes sense over any (computationally feasible)
ordered field.  Natural choices certainly include the field of rational numbers as well as certain (or even all)
algebraic extensions.\footnote{For a more general perspective on the subject see the book of Blum et
  al.~\cite{MR99a:68070}.}  In the following we focus on the rationals and their extensions by radical expressions.

Examining the Algorithms~\ref{alg:extract} and~\ref{alg:extend} more closely reveals that computations within the
coordinate domain are necessary only to decide whether a given point is contained in a given affine hyperplane and, if
not, on which side it lies.  That is to say, if we are only interested in combinatorially correct output, it suffices to
evaluate signs of determinants without ever knowing their precise values.  This approach, sometimes called \emph{robust}
or \emph{exact geometric computation}, is taken in several computational geometry libraries including
\texttt{CGAL}~\cite{cgal} and \texttt{LEDA}~\cite{leda}.  It has the advantage that it can be extended to radical
extensions of~$\QQ$ rather easily.  This is implemented in \texttt{LEDA} and in the \texttt{CORE} library~\cite{core}.

If, on the other hand, we do need the exact facet normal vectors then we have to use an exact implementation of the
arithmetic functions of our coordinate domain.  Today's standard for long integer and rational arithmetic is the
\texttt{GNU Multiprecision Library (GMP)}~\cite{gmp}.  In fact, this is the arithmetic used in \texttt{polymake}'s
implementation of Beneath-and-Beyond as well as in the codes of Avis~\cite{lrs} and Fukuda~\cite{cdd}.  There are
packages in computer algebra and computational number theory which can perform computations in (arbitrary) finite
extensions of~$\QQ$, but currently this functionality does not seem to be available as a stand-alone library which
encapsulates a given field as a number type to be used in standard C or C++ code.

Even if one agrees to compute with rational coordinates, on a technical level there is still one choice to be made: One
can either compute with representations of rational numbers directly, or one can translate everything into integer
coordinates, essentially by scaling.  While this usually does not make much of a difference for the algorithms, the
actual numbers that occur during the computations are different and, in particular, of different sizes.  This may affect
the performance of an algorithm, but often one does not know in advance which method is superior to the other.

It would be interesting to combine the two techniques, exact geometric computation and exact coordinates, for the
Beneath-and-Beyond algorithm in the following way: Use exact geometric computation to produce the combinatorially
triangulation, and only in the end compute the facet normals by solving systems of linear equations defined by the
$(d-1)$-faces on the boundary.  Since in many cases the arithmetic consumes most of the running time it should be
possible to save some time this way.  The author is not aware of an implementation of any convex hull code taking this
approach.

\section{Empirical Results}
\label{sec:empirical}
The analysis in Section~\ref{sec:sizes} showed that, from a theoretical point of view, Beneath-and-Beyond seems to
be a weak algorithm: It has a super-polynomial total running time even on simplicial polytopes.  In this section we
display a few computational results, where the performance of \texttt{polymake's} implementation
\texttt{beneath\_{}beyond} is compared to two other programs: Fukuda's \texttt{cdd}, which implements (dual)
Fourier-Motzkin elimination, see Ziegler~\cite[Sections 1.2 and~1.3]{GMZ} or Fukuda and
Prodon~\cite{MR98c:90108}, and Avis'~\texttt{lrs}, which uses reverse search, see Avis and
Fukuda~\cite{752.68082,lrs-revised}.

\begin{figure}
  \begin{center}
    \includegraphics[height=12cm,angle=270]{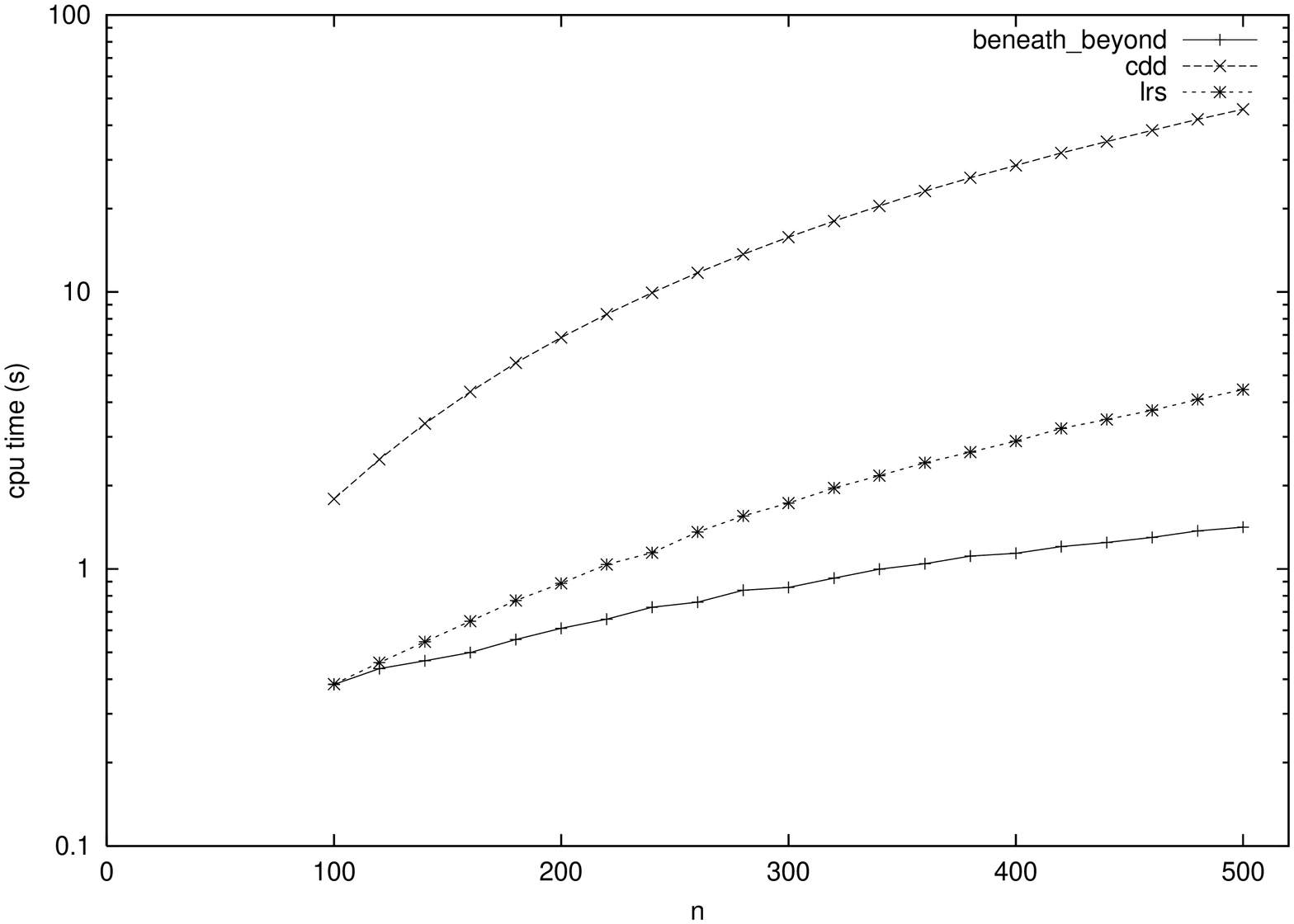}

    \bigskip
    \includegraphics[height=12cm,angle=270]{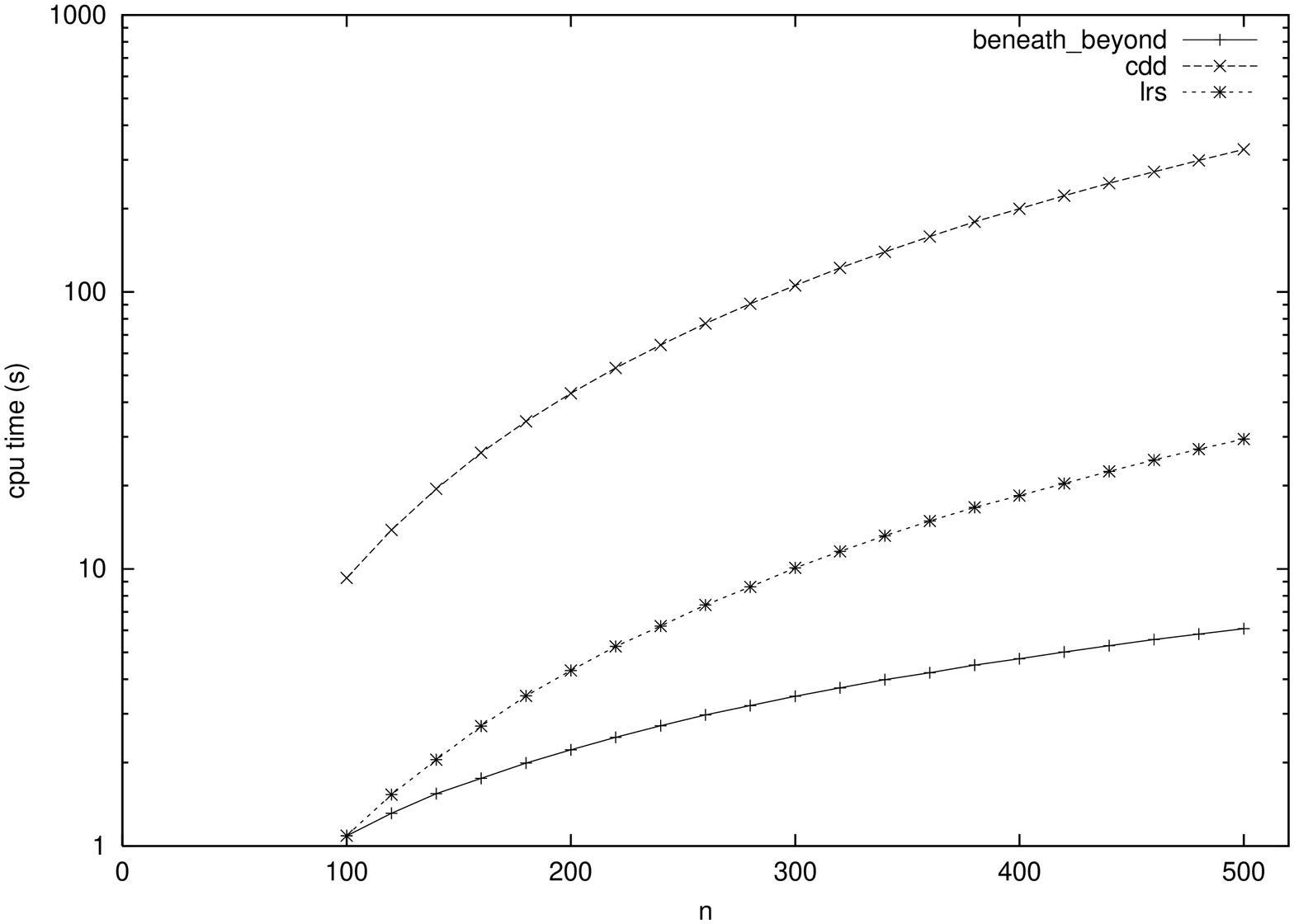}

    \caption{``Random spheres'' with $n$~vertices in dimensions~$3$ (top) and~$4$ (bottom).  Average over
      10 polytopes, each program run only once.  Note that our timings below 1 second are not very
      accurate.\label{fig:rand-sphere34}}
  \end{center}  
\end{figure}

\begin{figure}
  \begin{center}
    \includegraphics[height=12cm,angle=270]{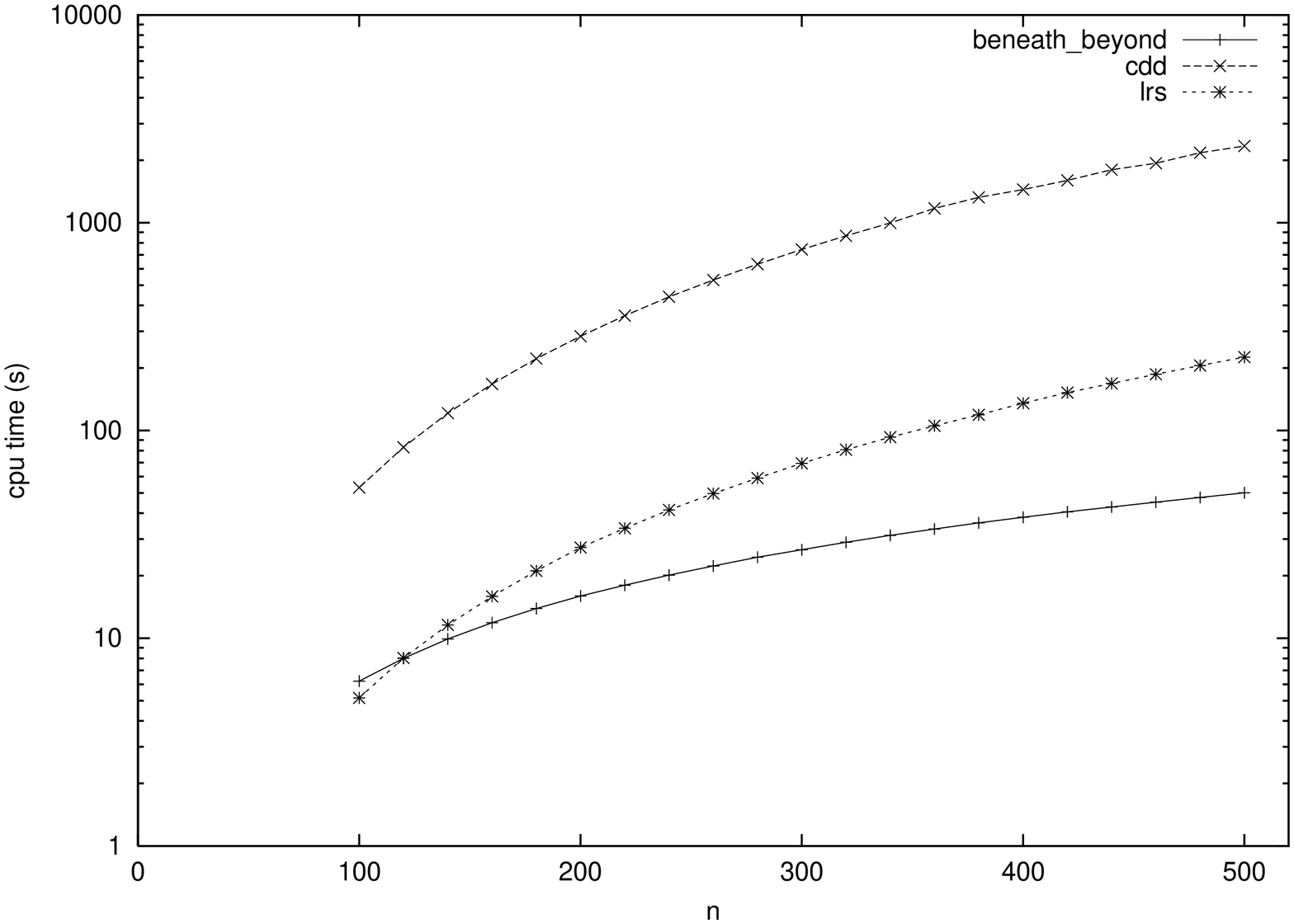}

    \bigskip
    \includegraphics[height=12cm,angle=270]{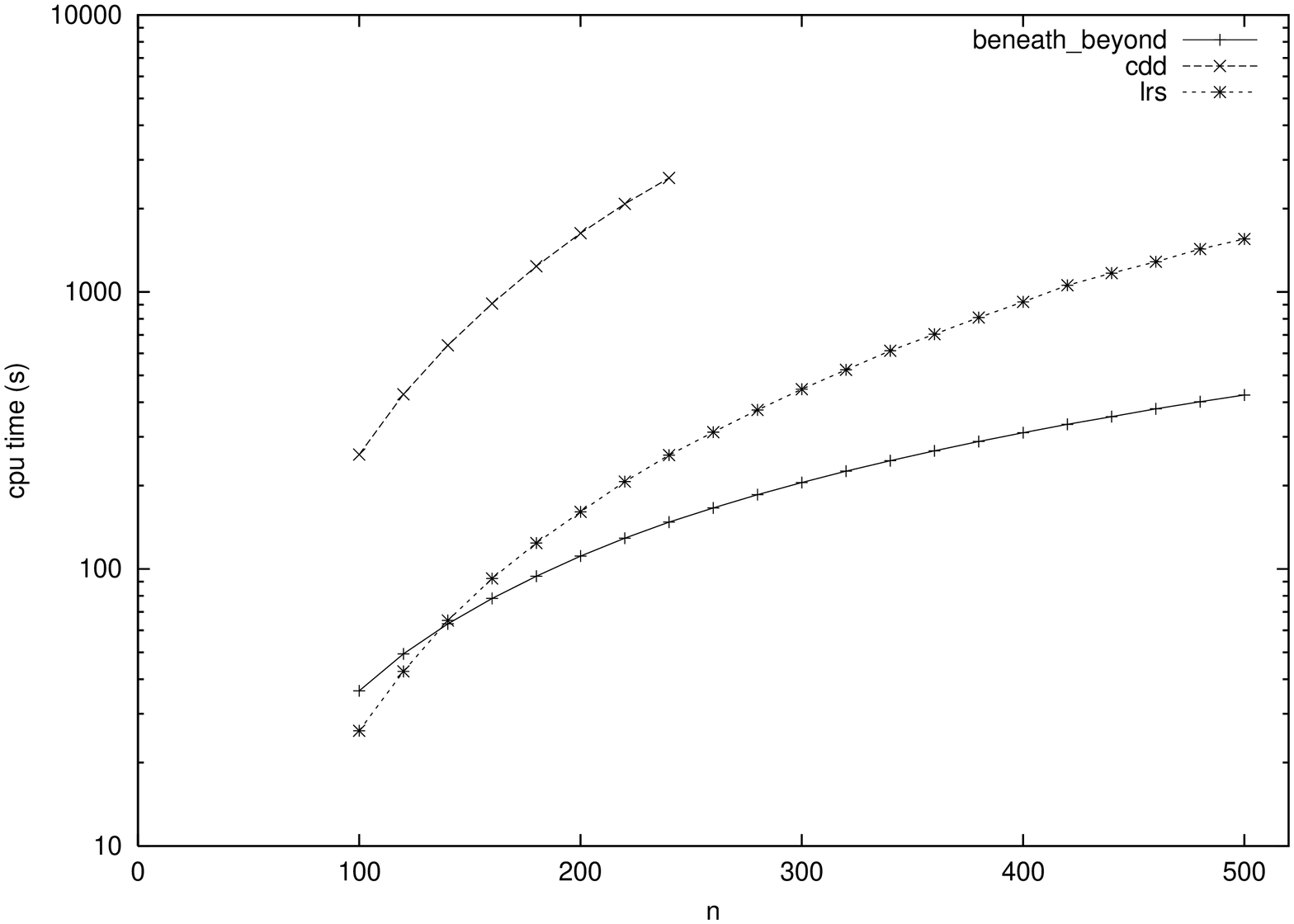}

    \caption{``Random spheres'' with $n$~vertices in dimensions~$5$ (top) and~$6$ (bottom).  Average over
      10 polytopes, each program run only once.  \texttt{cdd} not tested for input with more than $240$~vertices since
      it takes about three hours per test.\label{fig:rand-sphere56}}
  \end{center}  
\end{figure}

The computational results below are intended to complement the corresponding data in the paper by Avis, Bremner, and
Seidel~\cite{877.68119}.  The performance comparison among the three programs is fair in the sense that all programs use
the same implementation of exact rational arithmetic, namely from the \texttt{GNU Multiprecision Library
  (GMP)}~\cite{gmp}.  However, \texttt{beneath\_beyond} and \texttt{cdd} both use rational coordinates while
\texttt{lrs} uses integers only.  This could contribute to \texttt{beneath\_beyond}'s superiority over \texttt{lrs} for
the ``random spheres'', see Figures~\ref{fig:rand-sphere34} and~\ref{fig:rand-sphere56}; but this needs further
investigation.

Roughly speaking the input polytopes come in two groups: The first group consists of the intricate polytopes discussed
in the previous section.  These are the dwarfed cubes, products of simplices, and dwarfed products of polygons.  While
they are known to be computationally hard for all iterative convex hull algorithms, this study shows that
\texttt{beneath\_{}beyond} performs particularly bad.  The second class of polytopes investigated are convex hulls of
(uniformly distributed) random points on the unit sphere, sometimes called ``random spheres''.  Such polytopes are
almost always simplicial (in fact, all the tested ones were simplicial).  The performance of \texttt{beneath\_{}beyond}
is clearly better than that of \texttt{lrs}, which in turn clearly beats~\texttt{cdd}.

The ``random spheres'' were constructed with \texttt{polymake}'s client program \texttt{rand\_{}sphere} which produces
uniformly distributed random points on the unit sphere with double (i.e., IEEE 64-Bit floating point) coordinates
rounded to six digits after the decimal point, which were then transformed into exact rational numbers; the numerators
and denominators of such rational numbers typically have 15 to 20 decimal digits each.

While in~\cite{877.68119} the authors also compared various insertion strategies for the incremental algorithms we did
not do that here: \texttt{cdd} was run with the lexmin insertion rule (which is its default behavior).  This seems
justified in view of the computational results in loc.~cit., which showed that lexmin is optimal among the strategies
tested in most cases.

All timings are taken on several identical Linux machines with an Athlon XP 1800+ processor (1533.433 MHz clock, 3060.53
bogomips) and 512~MB main memory running RedHat Linux~7.3.  The timings were taken via Perl's \texttt{time()} function,
and we counted the time in user mode only.  The parameter \texttt{datasize} was limited to 400~MB.  Since all programs
were run through \texttt{polymake}, by calling the clients \texttt{cdd\_{}ch\_{}client} and \texttt{lrs\_{}ch\_{}client},
respectively, there is a certain additional overhead (due to socket communication and data conversion) which should
approximately be the same for the three programs.  This makes up for the slight difference in shape of some of the
curves as compared to the corresponding ones in~\cite{877.68119}.  In particular, timings below one second are almost
impossible to interpret this way.  Moreover, the implementations of \texttt{cdd} and \texttt{lrs} tested here are more
advanced than the ones tested in loc.~cit.  In particular, the old versions of \texttt{cdd} and \texttt{lrs} did not use
the (considerably faster) \texttt{GMP} arithmetic.

Due to the fact that the timings include some input/output overhead it is impossible to reproduce exact timing data.
Therefore we took the average over several runs. Additionally, we randomly permuted the order of the input.  While
\texttt{cdd} and \texttt{lrs} are almost insensitive to the ordering, this can make a considerable difference for
\texttt{beneath\_{}beyond}.  The highest variation among the individual timings of \texttt{beneath\_{}beyond} occurred
for the dwarfed products of polygons: A more detailed statistical analysis is given in
Table~\ref{tab:dwarfed-product-polygons}.

\begin{figure}
  \begin{center}
    \includegraphics[height=12cm,angle=270]{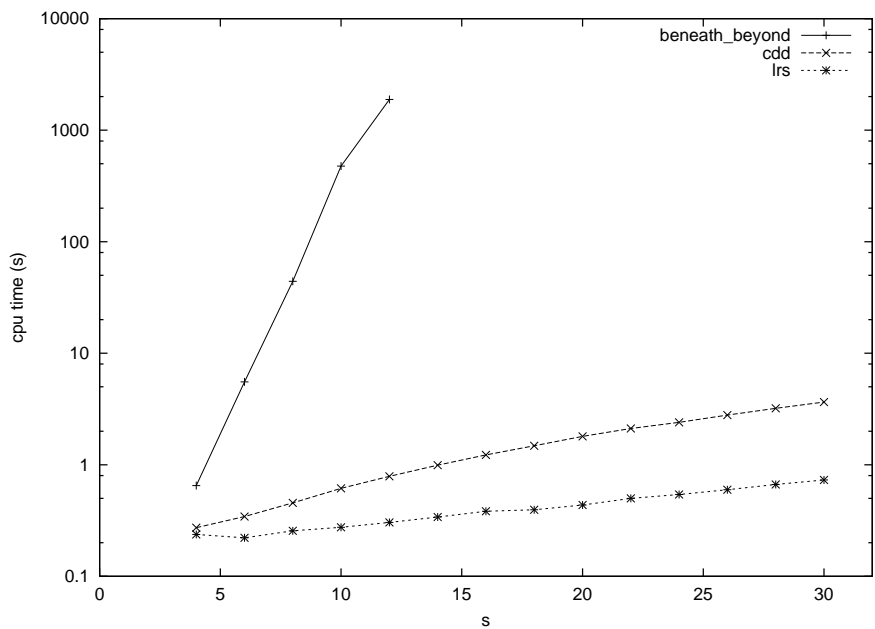}

    \caption{Dwarfed products $g_{10}(s)$ of five $s$-gons.  In this test $s$ is always even.
      Compare~\cite[Figure~5]{877.68119}.  \texttt{cdd} and \texttt{lrs}: average over 10 runs,
      \texttt{beneath\_beyond}: average over 50 runs.  Memory overflow (more than 400~MB required) in
      \texttt{beneath\_{}beyond} for $s>12$.  For more details of the \texttt{beneath\_{}beyond} timings see
      Table~\ref{tab:dwarfed-product-polygons}.\label{fig:dwarfed-product-polygons}}
  \end{center}  
\end{figure}

\begin{table}
  \caption{The timings (in seconds) for the individual runs of \texttt{beneath\_{}beyond} for the dwarfed products of
    polygons~$g_{10}(s)$, see Figure~\ref{fig:dwarfed-product-polygons}, vary rather strongly.  $50$ runs with a random
    insertion order were performed.\label{tab:dwarfed-product-polygons}}
  \begin{center}
    \footnotesize
    \begin{tabular*}{\textwidth}{@{\extracolsep{\fill}}rrrrr@{}}
      \toprule
      $s$ & average & minimum & maximum & standard deviation\\
      \midrule
      4 & 0.650 & 0.35 & 1.05 & 0.1628\\
      6 & 5.530 & 0.78 & 23.20 & 4.6607\\
      8 & 44.091 & 3.26 & 226.55 & 40.7324\\
      10 & 476.051 & 10.00 & 1984.59 & 452.4177\\
      12 & 1883.964 & 239.46 & 8560.20 &  1803.4397\\
      \bottomrule
    \end{tabular*}
  \end{center}
\end{table}

\begin{figure}
  \begin{center}
    \includegraphics[height=12cm,angle=270]{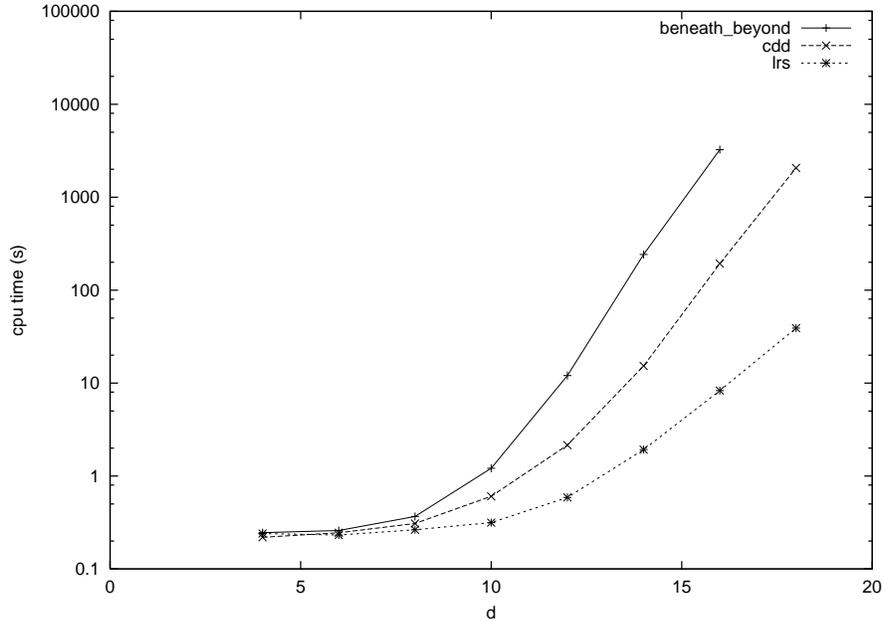}
    \caption{Dwarfed $d$-cubes~$c_d$, for $d$ even.  Compare~\cite[Figure~4]{877.68119}.  \texttt{cdd} and \texttt{lrs}: average over
      10 runs, \texttt{beneath\_beyond}: average over 50 runs, memory overflow (more than 400~MB required) for
      $d>16$.\label{fig:dwarfed-cubes}}
  \end{center}  
\end{figure}

\begin{figure}
  \begin{center}
    \includegraphics[height=12cm,angle=270]{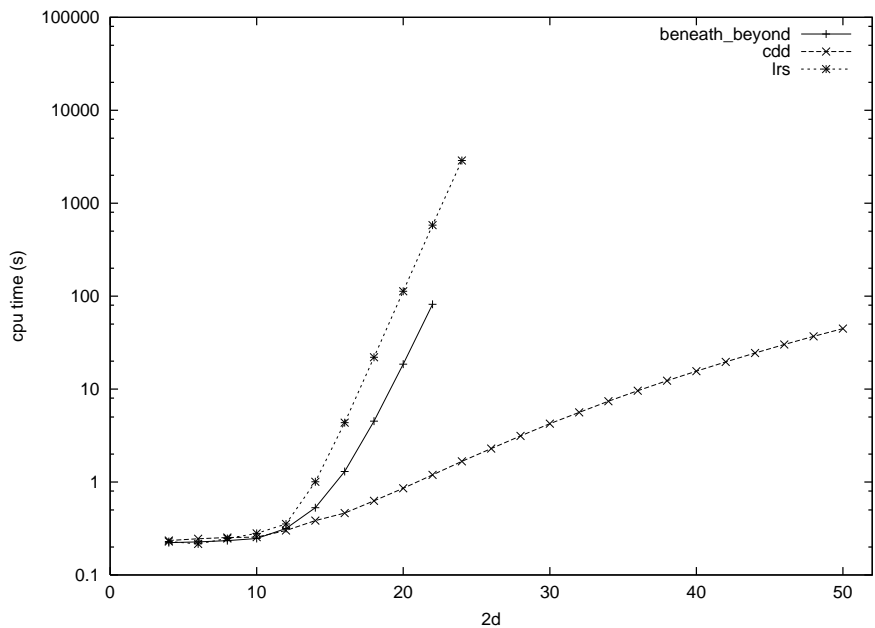}

    \caption{Products of two $d$-simplices.  Compare~\cite[Figure~2]{877.68119}.  Average over 10 runs for each program.
      \texttt{beneath\_{}beyond}: Memory overflow (more than 400~MB required) for $d>11$. \label{fig:products-simplices}}
  \end{center}  
\end{figure}

For the ``random spheres'' we changed our experimental set-up slightly.  For each set of parameters $d\in\{3,4,5,6\}$
and $s\in\{100,120,140,\dots,500\}$ we produced $10$ polytopes as the convex hulls of $s$ random points on the unit
sphere~$\Sph^{d-1}\subset\RR^d$.  Each of the three convex hull codes was run once on each polytope (for
\texttt{beneath\_{}beyond} we inserted the points in the same order as they were produced).  The charts show the average
values for each code taken over the $10$~samples.  The individual timings varied only by a little: The highest
deviations occurring for \texttt{cdd} and $d=6$, see the Table~\ref{tab:rand-sphere6-cdd} for more details.

\begin{table}
  \caption{The timings (in seconds) for the individual runs of \texttt{cdd} for $s$ random points on the unit sphere
    in~$\RR^6$, see Figure~\ref{fig:rand-sphere56}, vary in a rather narrow range.  We checked $10$ cases for each value
    of~$s$, produced independently at random.\label{tab:rand-sphere6-cdd}}
  \begin{center}
    \footnotesize
    \begin{tabular*}{\textwidth}{@{\extracolsep{\fill}}rrrrr@{}}
      \toprule
      $s$ & average & minimum & maximum & standard deviation\\
      \midrule
      100 & 258.997 & 247.05 & 273.58 & 7.5541\\
      120 & 427.075 & 416.67 & 446.39 & 8.3978\\
      140 & 641.947 & 627.80 & 666.61 & 12.0566\\
      160 & 907.759 & 895.52 & 938.18 & 12.9157\\
      180 & 1238.253 & 1202.69 & 1265.81 & 18.9353\\
      200 & 1630.816 & 1591.89 & 1670.96 & 23.6467\\
      220 & 2078.458 & 2028.65 & 2132.95 & 29.9611\\
      240 & 2582.256 & 2495.02 & 2650.31 & 41.1209\\
      \bottomrule
    \end{tabular*}
  \end{center}
\end{table}

\section{Concluding Remarks}
Authors in discrete geometry often talk about ``triangulations of finite point sets''.  They refer to triangulations of
the convex hull of these points.  For more details see Pfeifle and Rambau~\cite{AGSS:PfeifleRambau}.

Finite simplicial complexes, as combinatorial abstractions of triangulations, are an indispensible tool in topology.  At
first sight it may seem accidental that notions from topology appear in questions concerning the complexity of convex
hull computations.  However, as developed in the paper~\cite{OraclesConvexHulls} there is a deeper connection: Up to
polynomial equivalence a convex hull computation can be replaced by suitable simplicial homology computations.  For
details on algorithms and implementations to compute homology see Dumas et al.~\cite{AGSS:DumasEtAl}.  The
paper~\cite{AGSS:KaibelPfetsch} by Kaibel and Pfetsch contains more information about the complexity status of the
convex hull problem.


\end{document}